
\def\LaTeX{%
  \let\Begin\begin
  \let\End\end
  \let\salta\relax
  \let\finqui\relax
  \let\futuro\relax}

\def\UK{\def\our{our}\let\sz s}
\def\USA{\def\our{or}\let\sz z}

\UK



\LaTeX

\USA


\salta

\documentclass[twoside,12pt]{article}
\setlength{\textheight}{24cm}
\setlength{\textwidth}{16cm}
\setlength{\oddsidemargin}{2mm}
\setlength{\evensidemargin}{2mm}
\setlength{\topmargin}{-15mm}
\parskip2mm


\usepackage[usenames,dvipsnames]{color}
\usepackage{amsmath}
\usepackage{amsthm}
\usepackage{amssymb, bbm}
\usepackage[mathcal]{euscript}
\usepackage{cite}

\usepackage{hyperref}
%
%


\definecolor{viola}{rgb}{0.3,0,0.7}
\definecolor{ciclamino}{rgb}{0.5,0,0.5}
\definecolor{rosso}{rgb}{0.85,0,0}

\def\betti #1{#1}
\def\pier #1{#1}


\bibliographystyle{plain}


%

\finqui

\def\Beq{\Begin{equation}}
\def\Eeq{\End{equation}}

\def\Bthm{\Begin{theorem}}
\def\Ethm{\End{theorem}}
\def\Blem{\Begin{lemma}}
\def\Elem{\End{lemma}}
\def\Bprop{\Begin{proposition}}
\def\Eprop{\End{proposition}}

\def\Brem{\Begin{remark}\rm}
\def\Erem{\End{remark}}

\def\Bdim{\Begin{proof}}
\def\Edim{\End{proof}}
\def\Bcenter{\Begin{center}}
\def\Ecenter{\End{center}}
\let\non\nonumber




\def\step #1 \par{\medskip\noindent{\bf #1.}\quad}


\def\Holder{H\"older}
\def\frechet{Fr\'echet}
\def\aand{\quad\hbox{and}\quad}

\def\lhs{left-hand side}
\def\rhs{right-hand side}


\def\multibold #1{\def\arg{#1}%
  \ifx\arg\pto \let\next\relax
  \else
  \def\next{\expandafter
    \def\csname #1#1#1\endcsname{{\boldsymbol #1}}%
    \multibold}%
  \fi \next}

\def\pto{.}

\def\multical #1{\def\arg{#1}%
  \ifx\arg\pto \let\next\relax
  \else
  \def\next{\expandafter
    \def\csname cal#1\endcsname{{\cal #1}}%
    \multical}%
  \fi \next}

\def\multigrass #1{\def\arg{#1}%
  \ifx\arg\pto \let\next\relax
  \else
  \def\next{\expandafter
    \def\csname gr#1\endcsname{{\mathbb #1}}%
    \multigrass}%
  \fi \next}


\def\multimathop #1 {\def\arg{#1}%
  \ifx\arg\pto \let\next\relax
  \else
  \def\next{\expandafter
    \def\csname #1\endcsname{\mathop{\rm #1}\nolimits}%
    \multimathop}%
  \fi \next}

\multibold
qwertyuiopasdfghjklzxcvbnmQWERTYUIOPASDFGHJKLZXCVBNM.

\multical
QWERTYUIOPASDFGHJKLZXCVBNM.

\multigrass
QWERTYUIOPASDFGHJKLZXCVBNM.

\multimathop
diag dist div dom mean meas sign supp .


\def\accorpa #1#2{\eqref{#1}--\eqref{#2}}
\def\Accorpa #1#2 #3 {\gdef #1{\eqref{#2}--\eqref{#3}}%
  \wlog{}\wlog{\string #1 -> #2 - #3}\wlog{}}


\def\separa{\noalign{\allowbreak}}

\def\graffe #1{\mathopen\{#1\mathclose\}}

\def\<#1>{\mathopen\langle #1\mathclose\rangle}
\def\norma #1{\mathopen \| #1\mathclose \|}

\def\[#1]{\mathopen\langle\!\langle #1\mathclose\rangle\!\rangle}

\def\iot {\int_0^t}
\def\ioT {\int_0^T}
\def\intQt{\int_{Q_t}}
\def\intQ{\int_Q}
\def\iO{\int_\Omega}

\def\dt{\partial_t}
\def\dn{\partial_\nnn}

\def\cpto{\,\cdot\,}

\def\checkmmode #1{\relax\ifmmode\hbox{#1}\else{#1}\fi}

\def\aeQ{\checkmmode{a.e.\ in~$Q$}}

\def\aet{\checkmmode{a.e.\ in~$(0,T)$}}

\def\aat{\checkmmode{for a.a.~$t\in(0,T)$}}


\def\erre{{\mathbb{R}}}




\def\genspazio #1#2#3#4#5{#1^{#2}(#5,#4;#3)}
\def\spazio #1#2#3{\genspazio {#1}{#2}{#3}T0}

\def\L {\spazio L}
\def\H {\spazio H}

\def\C #1#2{C^{#1}([0,T];#2)}


\def\Lx #1{L^{#1}(\Omega)}
\def\Hx #1{H^{#1}(\Omega)}

\def\LQ #1{L^{#1}(Q)}

\def\Ldue{\Lx 2}
\def\Linfty{\Lx\infty}

\def\Huno{\Hx 1}


\let\badphi\phi
\let\phi\varphi

\let\theta\vartheta

\let\eps\varepsilon
\let\Lam\Lambda

\let\TeXchi\chi                         
\newbox\chibox
\setbox0 \hbox{\mathsurround0pt $\TeXchi$}
\setbox\chibox \hbox{\raise\dp0 \box 0 }
\def\chi{\copy\chibox}



\let\emb\hookrightarrow
\def\CO{C_\Omega}
\def\cd{c_\delta}

\def\phie{\badphi_e}
\def\phir{\badphi_r}
\def\ui{u_i}
\def\ue{u_e}
\def\uie{(\ui,\ue)}
\def\kas{\kappa_s}
\def\kae{\kappa_e}
\def\kai{\kappa_i}
\def\kar{\kappa_r}
\def\kamin{\kappa_*}
\def\kamax{\kappa^*}
\def\uimax{\ui^{max}}
\def\uemax{\ue^{max}}

\def\sz{s_0}
\def\ez{e_0}
\def\iz{i_0}
\def\rz{r_0}
\def\uz{u_0}

\def\soluz{(s,e,i,r)}

\def\dtau{\delta^\tau\!}
\def\stau{s_\tau}
\def\etau{e_\tau}
\def\itau{i_\tau}
\def\rtau{r_\tau}
\def\htau{h_\tau}
\def\soluztau{(\stau,\etau,\itau,\rtau)}

\def\y #1#2{#1^{(#2)}} 
\def\sol #1{(\y s#1,\y e#1,\y i#1,\y r#1)} 

\def\soluzl{(\xi,\eta,\iota,\rho)}
\def\hi{h_i}
\def\he{h_e}
\def\hie{(\hi,\he)}
\def\sh{s_h} 
\def\eh{e_h}
\def\ih{i_h}
\def\rh{r_h}
\def\xih{\xi^h} 
\def\etah{\eta^h}
\def\iotah{\iota^h}
\def\rhoh{\rho^h}
\def\soluzh{(\xih,\etah,\iotah,\rhoh)}
\def\Phih{\Phi^h}

\def\uistar{\ui^*} 
\def\uestar{\ue^*}
\def\uiestar{(\uistar,\uestar)}
\def\sstar{s^*}
\def\estar{e^*}
\def\istar{i^*}
\def\rstar{r^*}
\def\soluzstar{(\sstar,\estar,\istar,\rstar)}

\def\uin{u_{i,n}}
\def\uen{u_{e,n}}
\def\sn{s_n}
\def\en{e_n}
\def\In{i_n}
\def\rn{r_n}
\def\soluzn{(\sn,\en,\In,\rn)}

\def\soluza{(p,q,w,z)} 

\def\Uad{\calU_{ad}}
\def\Upiu{\calU^+}

\def\Vp{V^*}

\def\normaV #1{\norma{#1}_V}

\Begin{document}


%
\title{Optimal control of a reaction-diffusion\\ model related \pier{to} the spread of COVID-19}
\author{}
\date{}
\maketitle
\Bcenter
\vskip-1cm
{\large\sc Pierluigi Colli$^{(1)}$}\\
{\normalsize e-mail: {\tt pierluigi.colli@unipv.it}}\\[.25cm]
{\large\sc Gianni Gilardi$^{(1)}$}\\
{\normalsize e-mail: {\tt gianni.gilardi@unipv.it}}\\[.25cm]
{\large\sc Gabriela Marinoschi$^{(2)}$}\\
{\normalsize e-mail: {\tt gabriela.marinoschi@acad.ro}}\\[.25cm]
{\large\sc Elisabetta Rocca$^{(1)}$}\\
{\normalsize e-mail: {\tt elisabetta.rocca@unipv.it}}\\[.45cm]
$^{(1)}$
{\small Dipartimento di Matematica ``F. Casorati'', Universit\`a di Pavia}\\
{\small and \pier{Research Associate at the} IMATI -- C.N.R. Pavia}\\ 
{\small via Ferrata 5, 27100 Pavia, Italy}\\[.2cm]
$^{(2)}$
{\small ``Gheorghe Mihoc-Caius Iacob'' Institute of Mathematical Statistics\\
and Applied Mathematics of the Romanian Academy}\\
{\small Calea 13 Septembrie 13, 050711 Bucharest, Romania}\\[.2cm]
\Ecenter

\Begin{abstract}
\noindent
\pier{This paper is concerned with the well-posedness and optimal control problem
of a reaction-diffusion system for an epidemic
Susceptible-Infected-Recovered-Susceptible (SIRS) mathematical model in
which the dynamics develops in a spatially heterogeneous environment. Using
as control variables the transmission rates $u_{i}$ and $u_{e}$ of contagion
resulting from the contact with both asymptomatic and symptomatic persons,
respectively, we optimize the number of exposed and infected individuals at
a final time $T$ of the controlled evolution of the system. More precisely,
we search for the optimal $u_{i}$ and $u_{e}$ such that the number of
infected plus exposed does not exceed at the final time a threshold value $%
\Lambda $, fixed {a} priori. We prove here the existence of optimal
controls in a proper functional framework and we derive the first-order
necessary optimality conditions in terms of the adjoint variables.}

\End{abstract}
\salta
\pagestyle{myheadings}
\newcommand\testopari{\sc Colli \ --- \ Gilardi \ --- \ Marinoschi \ --- \ Rocca}
\newcommand\testodispari{{\sc \pier{Optimal} control of a model for COVID-19}}
\markboth{\testopari}{\testodispari}
\finqui
%

\section{Introduction}
\label{Intro}
\setcounter{equation}{0}


\pier{Since 2020, the literature of epidemic mathematical models has been enriched with many papers proposing interesting modeling ideas} (see, e.g., \cite{Gatto,Giordano,Jha,Linka, Wang,Zohdi}) 
\pier{often based on compartmental models,
where the considered population is divided into \textquotedblleft
compartments\textquotedblright\ based on their qualitative characteristics
(like, e.g., \textquotedblleft susceptible\textquotedblright ,
\textquotedblleft infected\textquotedblright , \textquotedblleft
recovered\textquotedblright ), with different assumptions about the nature
and rate of transfer across compartments. These models have naturally allow
for the introduction of diffusion terms. For a recent overview of
mathematical models for virus pandemic we refer to \cite{Bellomo4}. Let us
note, however, that in most of the cases the models present in the
literature are based on systems of ordinary differential equations (ODEs) in
time, while here we explore a compartmental model based on partial
differential equations (PDEs), to accurately account for spatial variations
at the continuum level. We also mention that epidemic models including
spatial diffusion have been proposed and investigated since long time ago
(see, e.g., \cite{Mottoni-Orlandi-Tesei, Fitzgibbon-1,
Fitzgibbon-2, Webb}). Very recently, a new epidemic diffusion model with
nonlinear transmission rates and diffusion coefficients was introduced and
tested in \cite{Vig1, Vig2}, while in \cite{acgrr} the authors
proved well-posedness results for an initial-boundary value problem associated to a
variant of the compartmental model for COVID-19 studied in~\cite{Vig1, Vig2}.}

\pier{In the present contribution, we consider an optimal control problem for the following PDE-based model}
\begin{align}
  & \dt s 
  + \ui s\hskip1pt i
  + \ue s\hskip1pt e 
  - \div (\kas \nabla s)
  = \gamma r
   && \hbox{in $Q$}
  \label{Iprima}
  \\
  & \dt e
  - \ui s\hskip1pt i
  - \ue s \hskip1pt e 
  + \sigma e
  + \phie e
  - \div (\kae \nabla e)
  = 0
  &&\hbox{in $Q$}
  \label{Iseconda}
  \\
  & \dt i 
  + \phir i 
  - \div (\kai \nabla i)
  = \sigma e
  &&\hbox{in $Q$}
  \label{Iterza}
  \\
  & \dt r 
  - \phir i
  - \phie e 
  - \div (\kar \nabla r)
  = - \gamma r
  && \hbox{in $Q$}
  \label{Iquarta}  
  \\
  & \dn s = \dn e = \dn i = \dn r = 0
  &&\hbox{on $\Sigma$}
  \label{Ibc}
  \\
  & \soluz(0) = (\sz,\ez,\iz,\rz)
  &&\hbox{in $\Omega$}
  \label{Icauchy}
\end{align}
\Accorpa\Ipbl Iprima Icauchy
in the space time cylinder $Q := \Omega\times(0,T)$,  
where $\Omega$ denotes from now on  a bounded open subset of \pier{$\erre^d$ with $d=3$ at most}, $\nnn$ and $\dn$ are the outward unit normal field on the boundary $\Gamma:=\partial\Omega$ and the corresponding directional derivative, respectively. Finally, $\Sigma:=\Omega\times (0,T)$ denotes the space-time boundary. \pier{We point out that the system \Ipbl\ can be considered a simplification of the model in} \cite{Vig1, Vig2}.

In \Ipbl\ \pier{the variables $s,$ $e,$ $i,$ $r$ represent the susceptible
population, the exposed population, the infected population, and the
recovered population, respectively. Compared to \cite{Vig1}, here we
neglect both the newborn and the deceased individuals after reaching the
maximum life, by assuming that these populations have small sizes and do not
contribute essentially to the epidemic transmission during the time period $(0,T)$. 
Thus, in our model the natality rate and the natural mortality rate
(which is not related to the disease) are zero.}


\pier{On the other hand, the essential spreading parameters are the transmission
rates $u_{i}$ and $u_{e}$ of contagion resulting from contact with an
asymptomatic or symptomatic person respectively. These will play the role of
the control variables in the sequel and will be allowed to depend on both
space and time in a bounded way.}

\pier{We assume that the diffusion coefficients $k_{s},$ $k_{e},$ $k_{i},$ $k_{r}$
are different, one from the other, and are allowed to depend on the space and time
variables still in a bounded way, thus describing a spatial
heterogeneous diffusion and a time depending dynamics of the spread. In
particular, they will be needed to be bounded away from zero.}

The function $\gamma$ on the right hand side of equations \eqref{Iprima} and \eqref{Iquarta} \pier{represents the rate at which 
immunity is lost and recovered individuals turn in the susceptible compartment. This rate is assumed time depending being naturally an increasing function.} 

\pier{Finally, $\sigma $ denotes a parameter corresponding to the latency (or incubation period), while $\phie$ and $\phir $ are the rates at which
exposed and infected individuals, respectively, become recovered.}

\pier{We mention that, as observed in COVID-19 epidemic, the exposed may also
spread the disease, with a transmission rate $u_{e}$ possibly much smaller
that $u_{i},$ so that this model can be viewed as a
susceptible-infected-recovered-susceptible (SIRS) model.}

\pier{Let us comment a little about the main differences which are here in
comparison with the model introduced in \cite{acgrr, Vig1, Vig2}. First of all,
in \cite{acgrr, Vig1, Vig2}, the terms $u_{i}s \hskip1pt{i}$ and $u_{e}s\hskip1pt{e}$ were
premultiplied by a function $A(n)$, where $n$ denotes the sum of all
compartment populations. However, the behavior of $A(n)$ was assumed to be
close to a constant for significant values of $n$. Furthermore, here we
generalize those terms by adding a possible dependence of the controls $%
u_{i}$, $u_{e}$ on the space and time variables. On the other hand, we
simplify the system by assuming that the diffusion coefficients $k_{s},$ $k_{e},$ $k_{i},$ $k_{r}$ (which can depend on space and time) are
independent of $n$. This assumption is the first step for approaching the
existence of a nonlinear model in a future paper. The interested reader can
see \cite{acgrr} for the case in which the diffusion coefficients are allowed to depend on $n$. 
but in \cite{acgrr} they were assumed to be equal in all the equations. The problem
corresponding to different coefficients  $k_{s},$ $k_{e},$ $k_{i},$ $k_{r}$ that depend also on $n$ is
still open and deserves a different analysis, as said before.}


\pier{Finally, we do not consider here the birth and mortality rate (denoted by $\bar \alpha$ and $\bar\mu$ in \cite{Vig1,Vig2}) for a modeling reason: 
as already mentioned, natural birth and death rates do not influence the
disease spreading (or influence it in a very small amount the total population) and so they can be neglected. 
%
%
Moreover, with respect to  \cite{Vig1,Vig2} we omit here the dynamics of the deceased population due to the disease
itself, by agreeing, for simplification, that the number of dead individuals due to the disease within the time interval $(0,T)$ is small with respect to
the number of population in the other compartments.}

In the present contribution, we are interested, beside the well-posedness of a suitable variational formulation of the state system \Ipbl, 
to the study of the following optimal control problem: 
\Beq
  \hbox{\sl Minimize the cost functional $\calJ$ on the set $\Uad$,}
  \label{ctrlI}
\Eeq
where we define the set $\Uad$ of the admissible controls and the cost functional $\calJ$ as~follows
\begin{align}
  & \Uad := \graffe{\uie \in \pier{L^\infty (Q)^2  :\ 0\leq\ui\leq\uimax\ \hbox{and}\ \ 0\leq{}}\ue\leq\uemax\ \aeQ}
  \label{defUadI}
  \\
  & \calJ\uie
   := \frac 12 \iO \bigl( \bigl( (e+i)(T) - \Lam \bigr)^+ \bigr)^2
  + \frac 12 \intQ \bigl( |\ui|^2 + |\ue|^2 \bigr)\pier{.}
  \label{costI}
\end{align}
\pier{Here,} $e$ and $i$ are the components of the solution $\soluz$ to \Ipbl\ corresponding to $\uie$\pier{;}   $\uimax,\uemax$ \ are nonnegative given bounded functions and $\Lambda$ is a given threshold value for the number of exposed plus infected \pier{individuals} at the final time $T$. The objective of the optimal control problem is indeed to keep the number of exposed plus infected rates below a certain given threshold value at the final time of the process. For other works where similar cost functionals have been considered, the reader can refer, e.g., to \cite{CGLMRR} where a similar approach has been used in order to keep the prostate index in a prostate tumor growth model under a certain threshold value. 
 
\pier{The motivation of considering $u_{i}$ and $u_{e}$ as controls is due to
their essential role in transmission. They relate to the infection mechanism
which depend on the way people mix and contact other people,  including in
their structure the number of contacts between susceptible and infected or
exposed individuals, respectively. It may depend on time because of seasonal
influence or due to the dynamics of the disease and also on space due to the
local spatial particularities. They can be used as control variables since
the contact reduction due to the prevention social measures applied within
the intervention period $(0,T)$ can determine the decrease of their values
up to optimal values that ensure reaching a minimum value of the objective
functional. Optimal control are of high interest in epidemics and at this
point we can refer the reader to a very few recent papers dealing with such
a subject \cite{Balderrama, Giordano, GM-AMO, GM-DCDS} where control
problems or identification of various coefficients in ODE models (without
diffusion) have been studied, and \cite{Medhaoui} for a reaction-diffusion
model control problem.} 
 
The main result of our contribution is related to the existence of solutions to Problem~\eqref{ctrlI} for which first-order necessary optimality conditions are derived in terms of the adjoint variables.  In order to achieve this result, we first prove well-posedness and boundedness of solutions for the state system \Ipbl, then we get a continuous dependence result in suitable functional spaces that allow us to prove the differentiability of the control-to-state mapping. Then, we derive first-order necessary optimality conditions for Problem~\eqref{ctrlI} in terms of the adjoint variables.

The paper is organized as follows. 
In the next section, we list our assumptions and notations
and state our results.
The proofs of Theorems~\ref{Wellposedness} and \ref{Contdep} regarding the well-posedness of the state system are given in Section~\ref{WP}, while 
Section~\ref{CONTROL} is devoted to the study of the optimal control problem. 


\section{Statement of the problem and results}
\label{STATEMENT}
\setcounter{equation}{0}

In this section, we state precise assumptions and present our results.
First of all, the set $\Omega\subset\erre^d$, \pier{with $d=1,2,3,$} 
is~assumed to be bounded, connected and smooth.
Next, if $X$ is a Banach space, then $\norma\cpto_X$ denotes its norm,
with the only exception of the space $H$ defined below
and the $L^p$ spaces ($1\leq p\leq\infty$) constructed on $\Omega$ and~$Q$, 
whose norms will be indicated by $\norma\cpto$ (i.e., without any subscript)
and by~$\norma\cpto_p\,$, respectively. 
Moreover, for simplicity, the same symbol used for some norm in $X$ 
will also stand for the norm in any power of~$X$.
Furthermore, the symbol $X^*$ denotes the dual space of~$X$.
We also introduce 
\Beq
  H := \Ldue
  \aand 
  V := \Huno
  \label{defspazi}
\Eeq
and we adopt the framework of the Hilbert triplet $(V,H,\Vp)$
obtained by identifying $H$ with a subspace of $\Vp$ in the usual way, namely, in order 
that $\<z,v>=\iO zv$ for every $z\in H$ and $v\in V$,
where $\<\cpto,\cpto>$ is the duality pairing between $\Vp$ and~$V$. 

\vskip 2mm

Now, we list the assumptions we postulate on the structure of the system:
\begin{align}
  & \kas\,,\,\kae\,,\,\kai\,,\,\kar\, \in \LQ\infty
  \aand
  \kamin \leq \kas\,,\,\kae\,,\,\kai\,,\,\kar \leq \kamax
  \quad \aeQ
  \non
  \\
  & \quad \hbox{for some positive constants $\kamin$ and $\kamax$}
  \label{hpk}
  \\
  & \hbox{$\phie$, $\phir$ and $\sigma$ are positive constants}
  \label{hpcoeff}
  \\
  & \gamma \in L^\infty(0,T) \quad \hbox{is nonnegative}.
  \label{hpgamma}
\end{align}
\Accorpa\HPstruttura hpk hpgamma

For the data, we make the assumptions listed below.
Even though the controls have to be considered as fixed data at the present stage,
it is convenient to introduce a constant $M$ which is an upper bound for their $L^\infty$ norms.
Thus, we require~that:
\begin{align}
  & \sz \,,\, \ez \,,\,\iz \,,\,\rz \in \Linfty
  \quad \hbox{are nonnegative}
  \label{hpz}
  \\
  & \ui \,,\, \ue \in \LQ\infty
  \quad \hbox{with} \quad 
  0 \leq \ui \,,\, \ue \leq M
  \quad \aeQ \,.
  \label{hpu}
\end{align}
\Accorpa\HPdati hpz hpu

The above assumptions guarantee that the problem \Ipbl\ 
has a unique variational solution
and that satisfactory stability and continuous dependence results hold true.

\Bthm
\label{Wellposedness}
Assume \HPstruttura\ on the structure of the system and \HPdati\ on the data.
Then there exists a unique quadruplet $\soluz$ enjoying the regularity properties
\begin{align}
  & s ,\, e ,\, i ,\, r \in \H1\Vp \cap \L2V \emb \C0H
  \label{regsoluz}
  \\
  & s ,\, e ,\, i ,\, r \geq 0
  \quad \aeQ
  \label{possoluz}
  \\
  & s ,\, e ,\, i ,\, r \in \LQ\infty
  \label{bddsoluz}
\end{align}
\Accorpa\Regsoluz regsoluz bddsoluz
and satisfying the variational equations
\begin{align}
  & \< \dt s , v >
  + \iO \bigl( \ui \, s i + \ue \, s e \bigr) \, v
  + \iO \kas \nabla s \cdot \nabla v
  - \iO \gamma r \, v
  = 0
  \label{prima}
  \\
  & \< \dt e , v >
  - \iO \bigl( \ui \, s i + \ue \, s e \bigr) \, v
  + \iO (\sigma+\phie) e \, v
  + \iO \kae \nabla e \cdot \nabla v
  = 0
  \label{seconda}
  \\
  & \< \dt i , v >
  + \iO \phir \, i \, v
  + \iO \kai \nabla i \cdot \nabla v
  - \iO \sigma e \, v
  = 0
  \label{terza}
  \\
  & \< \dt r , v >
  - \iO \bigl( \phir i + \phie e \bigr) \, v
  + \iO \kar \nabla r \cdot \nabla v
  + \iO \gamma r \, v
  = 0
  \label{quarta}
\end{align}
\aet\ and for every $v\in V$,
as well as the initial condition
\Beq
  \soluz(0) = (\sz,\ez,\iz,\rz).
  \label{cauchy}
\Eeq
\Accorpa\Pbl prima cauchy
Moreover, the stability estimate
\Beq
  \norma\soluz_{\H1\Vp\cap\C0H\cap\L2V\cap\LQ\infty}
  \leq K_1
  \label{stab}
\Eeq
holds true with some constant $K_1>0$ that depends only on the structure of the system,
$\Omega$, $T$, the initial data, and the constant~$M$.
In particular, $K_1$~is independent of~$\uie$.
\Ethm

\Bthm
\label{Contdep}
Under the assumptions of Theorem~\ref{Wellposedness} on the structure and the initial data,
let $\ui^{(j)},\,\ue^{(j)}$, $j=1,2$, be nonnegative functions in $\LQ\infty$ whose norms are bounded by~$M$, 
and let $\sol j$ be the corresponding solutions.
Then the inequality
\begin{align}
  & \norma{\sol1-\sol2}_{\C0H\cap\L2V}
  \non
  \\
  & \leq K_2 \, \norma{(\ui^{(1)},\ue^{(1)})-(\ui^{(2)},\ue^{(2)})}_{\L2H}
  \label{contdep}
\end{align}
holds true with a positive constant $K_2$ that depends only on the structure of the system,
$\Omega$, $T$, the initial data, and the constant~$M$.
\Ethm

By accounting for the above results, we can address the control problem sketched in the Introduction.
We set for convenience
\Beq
  \calU := (\LQ\infty)^2
  \aand
  \Upiu := \graffe{\uie\in\calU:\ \ui\geq0 \ \hbox{and} \ \ue\geq0 \ \aeQ}
  \label{defU}
\Eeq
and define the set $\Uad$ of the admissible controls and the cost functional $\calJ:\Upiu\to\erre$ as~follows
\begin{align}
  & \Uad := \graffe{\uie \in \Upiu:\ \ui\leq\uimax\ \hbox{and}\ \ \ue\leq\uemax\ \aeQ}
  \label{defUad}
  \\[2mm]
  & \calJ\uie
   := \frac 12 \iO \bigl( \bigl( (e+i)(T) - \Lam \bigr)^+ \bigr)^2
  + \frac 12 \intQ \bigl( |\ui|^2 + |\ue|^2 \bigr)
  \non
  \\
  & \quad \hbox{where $e$ and $i$ are the components of the solution $\soluz$}
  \non
  \\
  & \quad \hbox{to \Pbl\ corresponding to $\uie$}.
  \label{cost}
\end{align}
The quantities appearing in \eqref{defUad} and \eqref{cost} are required to satisfy
\begin{align}
  & \hbox{$\uimax,\uemax\in\LQ\infty$ \ are nonnegative}
  \label{hpUad}
  \\
  & \hbox{$\Lam$ \ is a positive constant}
  \label{hpLam}
\end{align}
and the symbol $(\cpto)^+$ denotes the positive part function. 
We notice that, for every $\uie\in\Upiu$, problem \Pbl\ actually has a unique solution by Theorem~\ref{Wellposedness},
so that $\calJ$ is well defined on~$\Upiu$.
Then, the control problem is the following:
\Beq
  \hbox{\sl Minimize the cost functional $\calJ$ on the set $\Uad$.}
  \label{ctrl}
\Eeq
In Section~\ref{CONTROL}, we prove the existence of an optimal control $(\uistar,\uestar)$
and derive first order optimality conditions.

\pier{\Brem
\label {Weights}
Of course, we could have put weights in the cost functional \pier{$\calJ$
by rescaling the integrals appearing in \eqref{cost} via suitable} positive coefficients.
In that case, \pier{the terms weigh differently according to their coefficients but it is still}
possible to develop the whole theory without any further difficulty.
\Erem}

The optimality condition we find sounds as follows:
if $(\uistar,\uestar)$ is an optimal control and $\soluzstar$ is the corresponding solution to \Pbl,
then the following variational inequality
\Beq
  \intQ \bigl( \sstar\istar(q-p) + \uistar \bigr) (\ui-\uistar)
  + \intQ \bigl( \sstar\estar(q-p) + \uestar \bigr) (\ui-\uistar)
  \geq 0
  \non
\Eeq
holds true for every $(\uistar,\uestar)\in\Uad$,
where $p$ and $q$ are the components of the weak solution $\soluza$ to a proper adjoint problem,
which is a variational backward parabolic problem whose strong form should be the following
\begin{align}
  & - \dt p - \div(\kas\nabla p) + (\uistar\istar+\uestar\estar) (p-q)
  = 0
  &&\hbox{in $Q$}
  \non
  \\
  & - \dt q - \div(\kae\nabla q) + (\sigma+\phie) q
    + \uestar\sstar (p-q) - \sigma w - \phie z
  = 0
  && \hbox{in $Q$}
  \non
  \\
  & - \dt w - \div(\kai\nabla w) + \phir (w-z)
    + \uistar\sstar(p-q) 
  = 0
  && \hbox{in $Q$}
  \non
  \\
  & - \dt z - \div(\kar\nabla z) + \gamma(z-p)
  = 0
  && \hbox{in $Q$}
  \non
  \\
  & \dn p = \dn q = \dn w = \dn z = 0
  && \hbox{on $\Sigma$}
  \non
  \\
  & p(T) = z(T) = 0
  \aand
  q(T) = w(T) = \bigl( (\estar+\istar)(T) - \Lam \bigr)^+
  && \hbox{in $\Omega$\,.}
  \non
\end{align}
However, this is formal and we do not refer to it in the sequel.
Indeed, our assumptions on the data do not ensure sufficient regularity, 
so that just the weak form of the adjoint problem can be considered.

\medskip

We conclude this section by recalling some tools 
and stating a general rule concerning the constants that appear in the estimates we perform in the sequel.
Throughout the paper, we will repeatedly use the Young inequality
\Beq
  a\,b \leq \delta\,a^2 + \frac 1{4\delta} \, b^2
  \quad \hbox{for all $a,b\in\erre$ and $\delta>0$}\,,
  \label{young}
\Eeq
as well as the \Holder\ and Schwarz inequalities.
Moreover, we take advantage of the three-dimensional continuous embedding
$V\emb\Lx p$ for $p\in[1,6]$, the embedding being compact if $p<6$,
and of the corresponding Sobolev and compactness inequalities
\begin{align}
  & \norma v_p
  \leq \CO \, \normaV v
  \quad \hbox{for every $v\in V$ and $p\in[1,6]$}
  \label{sobolev}
  \\
  & \norma v_p
  \leq \delta \norma{\nabla v} + C_{\Omega,\delta,p} \, \norma v
  \quad \hbox{for every $v\in V$, $p\in[1,6)$ and $\delta>0$}
  \label{compact}
\end{align}
where $\CO$ is a constant that depends only on~$\Omega$,
while $C_{\Omega,\delta,p}$ depends on $p$ and~$\delta$, in addition.

Finally, we will employ the abbreviations
\Beq
  Q_t := \Omega\times(0,t)
  \quad \hbox{for $t\in(0,T]$}
  \aand
  Q^t := \Omega\times(t,T)
  \quad \hbox{for $t\in[0,T)$}.
  \label{defQt}
\Eeq
Finally, here is our rule concerning the constants. 
The small-case symbol $c$ statnds for possibly different constants
(whose actual values may change from line to line and even within the same line)
that depend only on~$\Omega$, the structure of the system,
and the constants and the norms of the functions involved in the assumptions of the statements.
In particular, the values of $\,c\,$ may depend on the constant $M$ that appears in \eqref{hpu},
but they are independent of the control variables.
A~small-case symbol with a subscript like $\cd$
indicates that the constant may depend on the parameter~$\delta$, in addition.
On the contrary, we mark precise constants that we can refer~to
by using different symbols 
(e.g.,~a capital letter as in \eqref{sobolev}).


\section{The state system}
\label{WP}
\setcounter{equation}{0}

\pier{This section is dedicated to the study of the state system. We split the treatment in two steps,
where the former is concerned with the existence of solutions and uniform bound and the latter
deals with uniqueness and the continuous dependence estimate.}

\subsection{Existence}
\label{EXISTENCE}

This subsection regards the existence part of Theorem~\ref{Wellposedness}
and the stability estimate~\eqref{stab}.
Namely, we construct a solution $\soluz$ to \Pbl\ that satisfies~\eqref{stab}.
To this end, we first solve an approximating problem.
Given a positive integer~$n$, we set $\tau:=T/n$ and introduce the delay operator 
\begin{align}
  & \dtau : L^2(-\tau,T;H) \to \L2H
  \quad \hbox{defined by}
  \non
  \\
  & (\dtau v)(t) := v(t-\tau) 
  \quad \aat.
  \label{defdtau}
\end{align}
Then, the approximating problem consists in finding a quadruplet $\soluztau$
satisfying the analogue of \Regsoluz\
(with an obvious variant regarding~$\etau$, which must be defined even in $(-\tau,0)$)
and satisfying the variational equations
\begin{align}
  & \< \dt\stau , v >
  + \iO \bigl( \ui \, \stau \itau + \ue \, \stau \dtau\etau \bigr) \, v
  + \iO \kas \nabla\stau \cdot \nabla v
  - \iO \gamma \rtau \, v
  = 0
  \label{primatau}
  \\
  & \< \dt\etau , v >
  - \iO \bigl( \ui \, \stau \itau + \ue \, \stau \dtau\etau \bigr) \, v
  + \iO (\sigma+\phie) \etau \, v
  + \iO \kae \nabla\etau \cdot \nabla v
  = 0
  \label{secondatau}
  \\
  & \< \dt\itau , v >
  + \iO \phir \, \itau \, v
  + \iO \kai \nabla\itau \cdot \nabla v
  - \iO \sigma \dtau\etau \, v
  = 0
  \label{terzatau}
  \\
  & \< \dt\rtau , v >
  - \iO \bigl( \phir \itau + \phie \dtau\etau \bigr) \, v
  + \iO \kar \nabla\rtau \cdot \nabla v
  + \iO \gamma \rtau \, v
  = 0
  \label{quartatau}
\end{align}
\aet\ and for every $v\in V$,
as well as the initial condition
\Beq
  (\stau,\itau,\rtau)(0) = (\sz,\iz,\rz)
  \aand
  \etau(t) = \ez
  \quad \hbox{for $t\in[-\tau,0]$}.
  \label{cauchytau}
\Eeq
\Accorpa\Pbltau primatau cauchytau

We first recall the form of the weak maximum principle we use.
It is understood that an element $f\in\L2\Vp$ is nonnegative if
$\ioT\<f(t),v(t)>\,dt\geq0$ for every nonnegative $v\in\L2V$
and that, for $f,g\in\L2\Vp$, the inequality $g\leq f$ means that $f-g$ is nonnegative.

\Blem
\label{WeakMP}
Assume that $\kappa\in\LQ\infty$ is nonnegative, $\badphi\in\erre$,
$f\in\L2\Vp$ and $\uz\in H$ are both nonnegative,
and let $u\in\H1\Vp\cap\L2V$ be a solution to the problem
\begin{align}
  & \< \dt u , v > + \iO \badphi uv + \iO \kappa \nabla u \cdot \nabla v = \< f , v >
  \quad \hbox{\aet, for every $v\in V$}
  \non
  \\
  & u(0) = \uz \,.
  \non
\end{align}
Then $u\geq0$ \aeQ.
\Elem

We derive the following lemma

\Blem
\label{DaLSU}
Let $\kappa\in\LQ\infty$ satisfy $\kamin\leq\kappa\leq\kamax$ \aeQ\ and,
for given $f,g\in\L2\Vp$ and $\uz\in\Linfty$, let $u,w\in\H1\Vp\cap\L2V$ satisfy
\begin{align}
  & \< \dt u , v > + \iO \kappa \nabla u \cdot \nabla v = \< f , v >
  \quad \hbox{\aet, for every $v\in V$}
  \label{parab1}
  \\
  & \< \dt w , v > + \iO \kappa \nabla w \cdot \nabla v = \< g , v >
  \quad \hbox{\aet, for every $v\in V$}
  \label{parab2}
  \\
  & u(0) = w(0) = \uz .
  \label{pier1}
\end{align}
Then, we have that:
$i)$~If $f\in\L\infty H$ then $u\in\LQ\infty$ 
and its $L^\infty$ norm is bounded by a constant
depending only on $\Omega$, $T$, the constant $\kamin$ and~$\kamax$,
and the norms of the data $f$ and $\uz$.
$ii)$~If $f\in\L\infty H$, $g\leq f$ and $w\geq0$ \aeQ,
the same conclusions hold for~$w$.
\Elem

\Bdim
The first statement can be proved with minor changes
by following the argument used in \cite{LSU} to establish \cite[Thm.~7.1, p.~181]{LSU},
which regards the case of Dirichlet boundary conditions.
As for the second one, the inequality $g\leq f$ and the weak maximum principle 
imply that $w\leq u$. 
Since $w\geq0$ by assumption and $u$ is bounded by~$i)$, we deduce that $w$ is bounded as well
and that $\norma w_\infty\leq\norma u_\infty$,
whence the estimate of $\norma w_\infty$ \pier{follows} as desired.
\Edim

\Bthm
Problem \Pbltau\ has a unique solution $\soluztau$.
\Ethm

\Bdim
We solve the problem on the intervals $I_k:=[k\tau,(k+1)\tau]$ for $0\leq k<n$.
For simplicity, we still use the notation $\soluztau$ for each quadruplet that is defined just in~$I_k$.
For $k=0$ we recall that $\etau$ is defined on $[-\tau,0]$ by the second condition in~\eqref{cauchytau}.
So, on account of the whole \eqref{cauchytau},
we can construct a solution $\soluztau$ to the system on the interval $[0,\tau]$ 
by solving the linear parabolic variational equations 
\eqref{terzatau}, \eqref{quartatau}, \eqref{primatau} and \eqref{secondatau}
(and thus finding $\itau$, $\rtau$, $\stau$ and~$\etau$), in this order.
Notice that each component belongs to $H^1(0,\tau;\Vp)\cap L^2(0,\tau;V)$, thus to $C^0([0,\tau];H)$,
and is nonnegative by the weak maximum principle and bounded by \pier{Lemma~\ref{DaLSU}}
(applied with $\tau$ in place of~$T$),
since each equation can be presented in the form \eqref{parab1} with $f\in L^\infty(0,\tau;H)$.
However, no estimate of the $L^\infty$ norms is needed at the present level.
Once the problem is solved on~$[0,\tau]$, it is clear that the solution on the other intervals $[k\tau,(k+1)\tau]$
with $k<n$ can be constructed by induction on~$k$ by still considering
\eqref{terzatau}, \eqref{quartatau}, \eqref{primatau} and \eqref{secondatau} in this order.
Indeed, at each step, $\dtau\etau$~is a known function belonging to $L^\infty(k\tau,(k+1)\tau;H)$
and the initial value at $t=k\tau$ can be taken as the final value at $k\tau$ of the solution constructed in the previous step,
whose components belong to $C^0([(k-1)\tau,k\tau];H)$. 
So, the quadruplet $\soluztau$ constructed at the $k^{th}$ step belongs to
$(H^1(k\tau,(k+1)\tau;\Vp)\cap L^2(k\tau,(k+1)\tau;V))^4$, thus to $(C^0([k\tau,(k+1)\tau];H))^4$.
Moreover, each component is nonnegative and bounded,
once more for the weak maximum principle and  \pier{Lemma~\ref{DaLSU}} (where $(0,T)$ is replaced by the subinterval at hand).
This procedure constructs a quadruplet $\soluztau$ defined in the whole of~$[0,T]$
that is piecewise regular and global $H$-continuous
(due to the choice of the initial conditions at each step).
Moreover, it solves the system in each subinterval~$I_k$ and satisfies the initial conditions.
It follows that $\soluztau$ satisfies the global regularity we have required and is a global solution to \Pbltau.
Finally, uniqueness follows since each step provides a unique solution.
\Edim

Our aim is proving the existence of a solution $\soluz$ to \Pbl\ satisfying \eqref{stab}
by letting $\tau$ tend to zero (or~\pier{$n$} tend to infinity).
To this end, we need a number of a~priori estimates.

\step
First a priori estimate

We introduce the auxiliary function
\Beq
  \htau := \stau + \etau
  \label{defhtau}
\Eeq
and notice that it satisfies the variational equation
\begin{align}
  & \< \dt\htau , v >
  + \iO (\sigma + \phie) (\htau - \stau) v
  + \iO \bigl( \kas \nabla\stau + \kae \nabla(\htau-\stau) \bigr) \cdot \nabla v
  \non
  \\
  & = \iO \gamma \rtau \, v
  \quad \hbox{\aet, for every $v\in V$}.
  \label{auxil}
\end{align}
Now, we test \eqref{auxil} by $\htau$ and obtain that
\begin{align}
  & \frac 12 \, \frac d{dt} \iO |\htau|^2
  + \iO \kae |\nabla\htau|^2
  + \iO (\sigma + \phie) |\htau|^2
  \non
  \\
  & = - \iO (\kas-\kae) \nabla\stau \cdot \nabla\htau  \pier{{} + \iO (\sigma + \phie) \stau \htau}
  + \iO \gamma \rtau \htau
  \non
\end{align}
\aet.
At the same time, we test \eqref{primatau}, \eqref{terzatau} and \eqref{quartatau}
by $\lambda\stau$, $\itau$ and~$\rtau$, respectively,
where $\lambda$ is a positive parameter whose value will be chosen later~on.
Then, we add all the equalities we have obtained to each other and have~that
\begin{align}
  & \frac 12 \, \frac d{dt} \iO \bigl(
    |\htau|^2 + \lambda |\stau|^2 + |\itau|^2 + |\rtau|^2
  \bigr)
  \non
  \\
  & \quad {}
  + \iO \bigl(
     \kae |\nabla\htau|^2 + \lambda \kas |\nabla\stau|^2 + \kai |\nabla\itau|^2 + \kar |\nabla\rtau|^2 
  \bigr)
  \non
  \\
  & \quad {}
  + \iO (\sigma + \phie) |\htau|^2
  + \lambda \iO \bigl(
    \ui \itau |\stau|^2 + \ue \dtau\etau |\stau|^2
  \bigr)
  + \iO \phir |\itau|^2 
  + \iO \gamma |\rtau|^2
  \non
  \\
  & = - \iO (\kas - \kae) \nabla\stau \cdot \nabla\htau
    \pier{{} + \iO (\sigma + \phie) \stau \htau}
  + \iO \gamma \rtau \htau
  \non
  \\
  & \quad
   + \lambda \iO \gamma \rtau \stau
  + \sigma \iO \dtau\etau \, \itau
  + \iO \phir \itau \rtau 
  + \iO \phie \dtau\etau \, \rtau \,.
  \non
\end{align}
We recall assumptions \eqref{hpk} on the diffusion coefficients
and that all the components of the solution are nonnegative.
Hence, we only have to estimate the \rhs.
All the terms can be dealt with by the Young inequality and 
the first one can be estimated this~way
\begin{align}
  & - \iO (\kas - \kae) \nabla\stau \cdot \nabla\htau
  \leq \frac 12 \iO \kae |\nabla\htau|^2
  + \frac 12 \iO \frac{|\kas-\kae|^2}\kae \, |\nabla\stau|^2
  \non
  \\
  & \leq \frac 12 \iO \kae |\nabla\htau|^2
  + \frac 12 \iO \frac{4(\kamax)^2}\kamin \, |\nabla\stau|^2
  \leq \frac 12 \iO \kae |\nabla\htau|^2
  + \frac 12 \iO \lambda \kas |\nabla\stau|^2
  \non
\end{align}
the last inequality holding with the choice $\lambda=4(\kamax)^2/(\kamin)^2$.
At this point, we integrate over~$(0,t)$.
As for the terms involving $\dtau\etau$,
we observe~that (recall~\eqref{defQt})
\Beq
  \intQt |\dtau\etau|^2
  \leq \int_{-\tau}^t \Bigl(\iO |\etau(s)|^2 \Bigr) ds
  \leq \tau |\ez|^2 + 2 \intQt (|\htau|^2 + |\stau|^2).
  \non
\Eeq
Therefore, we can apply the Gronwall lemma and conclude~that
\Beq
  \norma{(\htau,\stau,\itau,\rtau)}_{\L\infty H\cap\L2V}
  \leq c 
  \quad \hbox{whence also} \quad
  \norma\etau_{\L\infty H\cap\L2V}
  \leq c \,.
  \label{primastima}
\Eeq

\step
Second a priori estimate

We recall the Sobolev embedding $V\emb\Lx6$ and \pier{the related  inequality \eqref{sobolev}. Therefore, it turns out that if $a\in \L2V$ and $b\in\L{\infty}H$ 
then the product $ab$ is in  $ \L2{\Lx{3/2}}$ with
 \begin{align}
  &\norma{ab}_{\L2{\Lx{3/2}}} = \left( \int_0^T \biggl| \iO |a|^{3/2}  | b|^{3/2} \biggr|^{4/3} \right)^{\!\!1/2} 
  \non
  \\
  & \leq \left( \int_0^T \biggl| \iO |a|^{6}  \biggr|^{1/3}      \biggl| \iO | b|^{2} \biggr| \right)^{\!\!1/2}
 \leq  \norma{a}_{\L2{\Lx6}} \norma{b}_{\L\infty{L^2(\Omega)}}
  \label{pier2}
\end{align}
 where the \Holder\ inequality with exponents $4$ and $4/3$ has been used in the first instance. Thus,  since $\norma{\dtau\etau}_{\L\infty H}    \leq   \max  \left\{  \norma{\ez} , \norma\etau_{\L\infty H}\right\}$ and  $\Lx{3/2} \emb \Vp $, 
we can recall} \eqref{primastima} to deal with the terms involving products that appear in~\eqref{primatau}. Namely, \pier{thanks to \eqref{pier2}} we have that
\Beq
  \ui\stau\itau \,,\, \ue\stau\dtau\etau \in \L2\Vp
  \non 
\Eeq
and there holds the estimate\pier{%
\begin{align}
  & \norma{\ui\stau\itau}_{\L2\Vp}
  + \norma{\ue\stau\dtau\etau}_{\L2\Vp}
  \non
  \\
  & \leq c\, \norma{\ui\stau\itau}_{\L2{\Lx{3/2}}}
  + c \,\norma{\ue\stau\dtau\etau}_{\L2{\Lx{3/2}}}
  \non
  \\
  & \leq c \, \norma\ui_\infty \, \norma\stau_{\L2{\Lx6}} \, \norma\itau_{\L\infty{L^2(\Omega)}}
  \non
  \\
  & \quad {}
  +c\, \norma\ue_\infty \, \norma\stau_{\L2{\Lx6}} \, \norma{\dtau\etau}_{\L\infty{L^2(\Omega)}}
  \leq c \,.
  \non
\end{align}}%
Hence, a bound in $\L2\Vp$ for $\dt\stau$ follows by comparison in~\eqref{primatau}.
The same argument applied to \eqref{secondatau} yields an estimate for~$\dt\etau$.
Finally, the analogous bounds for the other time derivatives follow directly by comparison in \accorpa{terzatau}{quartatau}.
We conclude~that
\Beq
  \norma{(\dt\stau,\dt\etau,\dt\itau,\dt\rtau)}_{\L2\Vp} \leq c \,.
  \label{secondastima}
\Eeq

\step
Third a priori estimate

We apply Lemma~\ref{DaLSU}, precisely~$ii)$, to \eqref{primatau} with 
$f=\gamma\rtau$ and $g=f-\ui\stau\itau-\ue\stau\dtau\etau$
by recalling that \pier{$\norma{f}_{\L\infty H}\leq c $} and $g\in\L2\Vp$.
We deduce~that
\Beq
  \stau \in \LQ\infty
  \aand
  \norma\stau_\infty \leq c \,.
  \non
\Eeq
Using this, we can apply the part $i)$ of Lemma~\ref{DaLSU} to \pier{\eqref{secondatau}} 
with $$f=\ui\stau\itau+\ue\stau\dtau\etau \pier{{}- (\sigma+\phie)\etau}$$ 
and deduce an analogous estimate.
By proceeding in a similar way for the other equations, we conclude~that
\Beq
  \norma\soluztau_\infty \leq c \,.
  \label{terzastima}
\Eeq

\step
Conclusion of the existence proof

On account of \accorpa{primastima}{terzastima}
and of \pier{well-known} weak, weak star and strong compactness results
(see, e.g., \cite[Sect.~8, Cor.~4]{Simon} for the latter),
we have~that
\begin{align}
  & \stau \to s , \quad
  \etau \to e , \quad
  \itau \to i 
  \aand
  \rtau \to r
  \non
  \\
  & \quad \hbox{weakly in $\H1\Vp\cap\L2V$ and weakly star in $\LQ\infty$,\quad whence}
  \non
  \\
  & \quad \hbox{weakly in $\C0H$, strongly in $\LQ p$ for $1\leq p<+\infty$ and \aeQ}
  \label{tautozero}
\end{align}
as $\tau\searrow0$ (at~least for a subsequence $\tau_k\searrow0$).
In particular the quadruplet $\soluz$ satisfies \Regsoluz,
the initial conditions \eqref{cauchy}
and estimate~\eqref{stab} (by~the lower semicontinuity of the norms).
We prove that $\soluz$ solves the variational equations \accorpa{prima}{quarta} as well.
To this end, we first show that
\Beq
  \dtau\etau \to e
  \quad \hbox{strongly in~$\L2H$}.
  \label{convdtau}
\Eeq
We extend $e$ to $[-\tau,T]$ by setting $e(t)=\ez$ for $t\in[-\tau,0]$
and have~that
\begin{align}
  & \norma{\dtau\etau-\pier{e}}_{\L2H}^2
  \leq 2 \, \norma{\dtau\etau-\dtau e}_{\L2H}^2
  + 2 \, \norma{\dtau e - e}_{\L2H}^2
  \non
  \\
  & = \pier{2}\int_0^{T-\tau} \norma{\etau(t)-e(t)}^2 \, dt
  + 2 \ioT \norma{\dtau e(t) - e(t)}^2 \, dt .
  \non
\end{align}
The integrals of the last line tend to zero,
the first one by \eqref{tautozero} and the second one since $e:[-\tau,T]\to H$ is uniformly continuous.
Thus \eqref{convdtau} is proved.
At this point, the convergence just established, our assumptions on $\uie$ and \eqref{tautozero} ensure~that
\Beq
  \ui \stau \itau \to \ui s i
  \aand
  \ue \stau \dtau\etau \to \ue s e 
  \quad \hbox{strongly in $\LQ p$ for $p\in[1,2)$}
  \non
\Eeq
By taking (for instance) $p=4/3$, we see that
we can let $\tau$ tend to zero in the integrated version of \eqref{primatau}
with arbitrary test functions $\pier{{}v \in{}} \L2V\cap\LQ4$
and obtain the integrated version of \eqref{prima} with the same test functions,
which is equivalent to \eqref{prima} itself.
The same argument can be applied to obtain \eqref{seconda}.
As for \eqref{terza} and \eqref{quarta} the situation is even simpler.
This concludes the proof of the existence of a solution satifying~\eqref{stab}.

\subsection{Uniqueness and continuous dependence}
\label{CONTDEP}

Here, we prove the uniqueness part of Theorem~\ref{Wellposedness} 
and the continuous dependence estimate presented in Theorem~\ref{Contdep}.
Our strategy is the following.
First, given $\ui^{(j)},\,\ue^{(j)}$, $j=1,2$, as in the statement,
we prove an auxiliary estimate similar to \eqref{contdep}
by assuming that $\sol j$ are arbitrary solutions corresponding to~$(\ui^{(j)},\,\ue^{(j)})$.
Then, we derive the uniqueness of the solution as a consequence.
Finally, we complete the proof of \eqref{contdep} with a constant $K_2$ with the dependence specified in the statement.

\step
An auxiliary estimate

We fix $\y\ui j,\,\y \ue j$, $j=1,2$, as in the statement
and assume that $\sol j$ are arbitrary solutions corresponding to~$(\y\ui j,\y\ue j)$.
We set for brevity
\begin{align}
  & \ui := \y\ui1 - \y\ui2 
  \aand
  \ue := \y\ue1 - \y\ue2 , \quad 
  \hbox{as well as}
  \non
  \\
  & s := \y s1 - \y s2 , \quad
  e := \y e1 - \y e2 , \quad
  i := \y i1 - \y i2 
  \aand
  r := \y r1 - \y r2 \,.
  \non
\end{align}
Now, we write each of the equations \accorpa{prima}{quarta} for both solutions and take the difference.
Then, we test the \pier{obtained} equalities by $s$, $e$, $i$ and~$r$, respectively, sum~up and rearrange.
We have~that
\begin{align}
  & \frac 12 \, \frac d{dt} \iO \bigl(
    |s|^2 + |e|^2 + |i|^2 + |r|^2
  \bigr)
  + \iO \bigl(
    \kas |\nabla s|^2 + \kae |\nabla e|^2 + \kai |\nabla i|^2 + \kar |\nabla r|^2
  \bigr)
  \non
  \\
  & \quad {}
  + \iO \bigl(
    (\sigma + \phie) |e|^2 + \phir |i|^2 + \gamma |r|^2
  \bigr)
  \non
  \\
  & = - \iO (\y\ui1 \y s1 \y i1 - \y\ui2 \y s2 \y i2) (s-e)
  - \iO (\y\ue1 \y s1 \pier{\y e1} - \y\ue2 \y s2 \pier{\y e2}) (s-e)
  \non
  \\
  & \quad {}
  + \iO (\gamma r s + \sigma e i + \phir i r + \phie e r).
  \non
\end{align}
At this point, we estimate from below the diffusion coefficients by $\kamin$ (recall~\eqref{hpk})
and estimate the terms on the \rhs.
However, we do it just for the first one, 
since the second one is analogous and the last one is simpler.
By first adding and subtracting suitable products and then estimating\pier{, with the help of the Young inequality}
we have~that
\begin{align}
  & - \iO (\y\ui1 \y s1 \y i1 - \y\ui2 \y s2 \y i2) (s-e)
  = - \iO (\y\ui1 s \y i1 \pier{{}+{}} \y\ui1 \y s2 i + \ui \y s2 \y i2) (s-e)
  \non
  \\
  & \leq \iO \bigl(
    M \, \norma{\y i1}_\infty \, |s|
    + M \, \norma{\y s2}_\infty \, |i| \pier{{}+ |\ui|\,\norma{\y s2}_\infty \norma{\y i2}_\infty}
  \bigr) (|s|+|e|)
  \non
  \\
  & \leq c \iO \bigl( |\ui|^2 + |s|^2 + |e|^2 + |i|^2 \bigr)
  \non
\end{align}
where the value of $c$ also depends on the $L^\infty$ norms of some of the components of the solutions.
By also treating the other integrals as said and then integrating over~$(0,t)$,
we see that we are in position to apply the Gronwall lemma.
By doing it, we conclude that
\Beq
  \norma\soluz_{\L\infty H\cap\L2V} 
  \leq c \, \norma{\uie}_{\L2H} 
  \label{auxiliary} 
\Eeq
where $c$ also depends on the solutions we have considered.

\step
Uniqueness

We apply \eqref{auxiliary} with $\y\ui1=\y\ui2$ and $\y\ue1=\y\ue2$
i.e., with $\ui=0$ and $\ue=0$, 
and obtain $\soluz=(0,0,0,0)$.
This means uniqueness, since the solutions $\sol1$ and $\sol2$ are arbitrary.

\step
Conclusion

We come back to the first step of the proof.
Because of uniqueness, the solution considered there must coincide 
with those constructed in the previous subsection.
Therefore, their norms we have considered 
are bounded by the constant $K_1$ of the statement of Theorem~\ref{Wellposedness}.
It follows that, a~posteriori, the constant $c$ appearing in \eqref{auxiliary}
can be replaced by a constant independent of the solutions,
so that \eqref{auxiliary} coincides with \eqref{contdep} and the proof is complete.


\section{The control problem}
\label{CONTROL}
\setcounter{equation}{0}

In this section, we address the control problem \eqref{ctrl}.
It is understood that all of the assumptions on the structure of the original system, the data,
and the ingredients of the cost functional~\eqref{cost}, which we have made throughout the paper, are in force from now~on.


\subsection{Existence of an optimal strategy}
\label{EXISTCONTROL}

The first result of ours is the following:

\Bthm
\label{OKcontrol}
The control problem \eqref{ctrl} admits at least one solution $(\uistar,\uestar)$.
\Ethm

\Bdim
We use the direct method.
Thus, let $\graffe{\uin,\uen}$ a minimizing sequence 
and $\graffe{\soluzn}$ the sequence of the corresponding states.
Since each $(\uin,\uen)$ belongs to~$\Uad$,
we can apply \eqref{stab} and deduce the corresponding bound for $\soluzn$.
Therefore, (at least for a not relabeled subsequence) we~have that
\begin{align}
  & \uin \to \ui
  \aand
  \uen \to \ue
  \quad \hbox{weakly star in $(\LQ\infty)^2$}
  \non
  \\
  & \sn \to s , \quad
  \en \to e , \quad
  \In \to i 
  \aand
  \rn \to r
  \non
  \\
  & \quad \hbox{weakly in $\H1\Vp\cap\L2V$ and weakly star in $\LQ\infty$,\quad whence}
  \non
  \\
  & \quad \hbox{weakly in $\C0H$, strongly in $\LQ p$ for $1\leq p<+\infty$ and \aeQ}.
  \non
\end{align}
Clearly, $\uie$ belongs to $\Uad$.
Moreover, by the same argument used to prove the existence of a solution
(see the last lines of Subsection~\ref{EXISTENCE}),
we see that $\soluz$ is the solution to \Pbl\ corresponding to~$\uie$.
Now, we \pier{show} a strong convergence property.
We write equation \accorpa{prima}{quarta} both for $\soluzn$ 
and for $\soluz$ with the data $(\uin,\uen)$ and $\uie$, respectively, and take the difference.
Then, we test the equality we obtain by $\sn-s$, $\en-e$, $\In-i$ and $\rn-r$, respectively.
Finally, we sum up and integrate over $(0,T)$.
We have~that
\begin{align}
  & \frac 12 \iO \bigl(
    |\sn(T)-s(T)|^2
    + |\en(T)-e(T)|^2
    + |\In(T)-i(T)|^2
    + |\rn(T)-r(T)|^2
  \bigr)
  \non
  \\
  & \quad {}
  + \intQ \bigl(
    \kas |\nabla(\sn-s)|^2
    + \kae |\nabla(\en-e)|^2
    + \kai |\nabla(\In-i)|^2
    + \kar |\nabla(\rn-r)|^2
  \bigr)
  \non
  \\
  & = \intQ (\uin \sn \In - \ui s i) (\en - e - \sn + s)
  + \intQ (\uen \sn \en - \ue s e) (\en - e - \sn + s)
  \non
  \\
  & \quad {}
  + \intQ (\gamma\rn - \gamma r) (\sn - s - \rn + r)
  - \intQ (\sigma + \phie) |\en - e|^2
  \non
  \\
  & \quad {}
  + \intQ [\sigma (\en - e) - \phir (\In - i)] (\In - i)
  + \intQ [\phir (\In - i) + \phie (\en - e)] (\rn - r) .
  \non
\end{align}
All the integrals on the \rhs\ tend to zero as $n$ tands to infinity since
\begin{align}
  & (\uen,\uin) \to (\ue,\ui)
  \quad \hbox{weakly star in $(\LQ\infty)^2$}
  \aand
  \non
  \\
  & \soluzn \to \soluz
  \quad \hbox{strongly in $(\LQ4)^4$}.
  \non
\end{align}
Therefore, the conclusions we are interested in is the following
\Beq
  \en(T) \to e(T) 
  \aand
  \In(T) \to i(T) 
  \quad \hbox{strongly in $H$}.
  \non
\Eeq
By also using the semicontinuity of the norm in $\LQ2$, we deduce that
\begin{align}
  & \liminf_{n\to\infty} \calJ(\en,\In)
  = \liminf_{n\to\infty} \Bigl\{
    \frac 12 \iO \bigl( \bigl( (\en+\In)(T) - \Lam \bigr)^+ \bigr)^2
    + \frac 12 \intQ \bigl( |\uin|^2 + |\uen|^2 \bigr)
  \Bigr\}
  \non
  \\
  & = \lim_{n\to\infty}\frac 12 \iO \bigl( \bigl( (\en+\In)(T) - \Lam \bigr)^+ \bigr)^2
  + \liminf_{n\to\infty} \intQ \bigl( |\uin|^2 + |\uen|^2 \bigr)
  \non
  \\
  & \geq \frac 12 \iO \bigl( \bigl( (e+i)(T) - \Lam \bigr)^+ \bigr)^2
  + \intQ \bigl( |\ui|^2 + |\ue|^2 \bigr)
  = \calJ(e,i).
  \non
\end{align}
Since $\graffe{(\uin,\uen)}$ is a minimizing sequence, we conclude that
$\calJ(e,i)$ is the minimum of $\calJ$ over~$\Uad$
and the proof is complete.
\Edim
 

\subsection{The control-to-state mapping}
\label{CTS}

In this section, we introduce the {\sl control-to-state} mapping $\calS$
and prove its \frechet\ differentiability in a suitable mathematical framework.
Namely, by recalling~\eqref{defU}, we~define
\begin{align}
  & \calY := \bigl( \H1\Vp \cap \L2V \bigr)^4
  \aand 
  \calS : \Upiu \to \calY
  \quad \hbox{by setting}
  \non
  \\& \hbox{$\calS\uie$ is the solution $\soluz$ to \Pbl\ corresponding to $\uie$}.
  \label{defS}
\end{align}
We notice that the domain of $\calS$ is not an open subset of~$\calU$.
Nevertheless, we can speak of \frechet\ differentiability by extending the usual definition to this case,
being understood that, for a given $\uie\in\Upiu$, 
the point $(\ui+\hi,\ue+\he)$ corresponding to the variation $(\hi,\he)\in\calU$ still belongs to~$\Upiu$.
Notice that this constraint does not compromise the uniqueness of the \frechet\ derivative,
since every $\hie\in\Upiu$ is an admissible variation.

Our differentiability result involves the linearized system we introduce at once.
Given any $\uie\in\Upiu$ and a variation $\hie\in\calU$,
and setting $\soluz=\calS\uie$,
the linearized system associated with $\uie$ and $\hie$ 
consists in finding a quadruplet $\soluzl\in\calY$ that satisfies the variational equations
\begin{align}
  & \< \dt\xi , v >
  + \iO \bigl(
    (\ui i)\xi 
    + (\ui s) \iota 
    + (\ue e)\xi 
    + (\ue s) \eta 
  \bigr) \, v
  \non
  \\
  & \quad {}
  + \iO \kas \nabla\xi \cdot \nabla v
  - \iO \gamma \rho \, v
  = - \iO \bigl( (si) \hi + (se) \he \bigr) \, v
  \label{primal}
  \\
  & \< \dt\eta , v >
  - \iO \bigl(
    (\ui i)\xi 
    + (\ui s) \iota 
    + (\ue e)\xi 
    + (\ue s) \eta 
  \bigr) \, v
  + \iO (\sigma + \phie) \eta \, v
  \non
  \\
  & \quad {}
  + \iO \kae \nabla\eta \cdot \nabla v
  = \iO \bigl( (si) \hi + (se) \he \bigr) \, v
  \label{secondal}
  \\
  & \< \dt \iota , v >
  + \iO \phir \iota \, v
  + \iO \kai \nabla\iota \cdot \nabla v
  - \iO \sigma \eta \, v
  = 0
  \label{terzal}
  \\
  & \< \dt\rho , v >
  - \iO \bigl( \phir \iota + \phie \eta \bigr) \, v 
  + \iO \kar \nabla\rho \cdot \nabla v
  + \iO \gamma \rho \, v
  = 0
  \label{quartal}
\end{align}
\aet\ and for every $v\in V$, 
as well as the initial condition
\Beq
  \soluzl(0) = (0,0,0,0) .
  \label{cauchyl}
\Eeq
\Accorpa\Pbll primal cauchyl
This is a uniformly parabolic problem system with bounded coefficients
since $\soluz$ is bounded by Theorem~\ref{Wellposedness}.
Hence, it has a unique solution $\soluzl$ which must satisfy
\Beq
  \norma\soluzl_\calY \leq c \, \norma\hie_{\pier{\L2H}}
  \non
\Eeq
with some constant~$c$.
However, for our purpouse, a weaker inequality is sufficient, 
as said in the next statement, where the precise dependence of the constant is specified.

\Bthm
\label{Wellposednessl}
For every $\uie\in\Upiu$ and $\hie\in\calU$,
the corresponding linearized system \Pbll\ has a unique solution $\soluzl\in\calY$.
Moreover, the following estimate
\Beq
  \norma\soluzl_\calY \leq C \, \norma\hie_\calU
  \label{perfrechet}
\Eeq
holds true with a positive constant $C$ that depends only on the structure of the original system,
$\Omega$, $T$, the initial data, and an upper bound for the norm~$\norma\uie_\infty$.
\Ethm

At this point, we are ready to prove the \frechet\ differentiability of~$\calS$.

\Bthm
\label{Frechet}
The control-to-state mapping $\calS$ is \frechet\ differentiable at every point of~$\Upiu$.
Namely, given any $\uie\in\Upiu$, the \frechet\ derivative $D\calS\uie$ is the operator
belonging to $\calL(\calU,\calY)$ that to every $\hie\in\calU$ associates the solution $\soluzl$ to the linearized system \Pbll\
corresponding to $\uie$ and the variation~$\hie$.
\Ethm

\Bdim
We fix $u:=\uie$ in $\Upiu$ once and for all
and first observe that the map $h:=\hie\mapsto\soluzl$ 
specified in the statement actually belongs to $\calL(\calU,\calY)$
thanks to Theorem~\ref{Wellposednessl}.
Thus, we only have to prove~that
\begin{align}
  & \norma{\calS(u+h) - \calS(u) - \soluzl}_\calY
  \leq \norma h_\calU \, \eps \bigl( \norma h_\calU \bigr)
  \non
  \\
  & \quad \hbox{whenever $u+h$ belongs to $\Upiu$, for some function}
  \non
  \\
  & \quad \hbox{$\eps:[0,+\infty)\to\erre$ that tends to zero at the origin}.
  \label{tesifrechet}
\end{align}
Without loss of generality, we can assume that $\norma h_\calU$ is small enough.
So, we fix $R>0$ larger than $\uimax$ and~$\uemax$ (see \eqref{defUad} and~\eqref{hpUad})
and assume that both $\norma{\ui+\hi}_\infty$ and $\norma{\ue+\he}_\infty$ do not exceed~$R$. 
Hence, we can apply \pier{Theorem~\ref{Wellposedness}} by choosing $M=R$ with both $u$ and $u+h$
and obtain the estimate \eqref{stab} for the corresponding solutions,~i.e.,
\begin{align}
  & \norma\soluz_{\H1\Vp\cap\C0H\cap\L2V\cap\LQ\infty}
  \leq K_1
  \aand
  \non
  \\
  & \norma{(\sh,\eh,\ih,\rh)}_{\H1\Vp\cap\C0H\cap\L2V\cap\LQ\infty}
  \leq K_1
  \non
\end{align}
where we have set 
\Beq
  \soluz := \calS(u) 
  \aand
  (\sh,\eh,\ih,\rh) := \calS(u+h) .
  \non
\Eeq
We also define
\Beq
  \soluzh := (\sh-s-\xi,\,,\eh-\pier{e}-\eta,\,,\ih-i-\iota,\,,\rh-r-\rho) 
  \non
\Eeq
and notice that $\soluzh$ belongs to $\calY$.
Moreover, we set for brevity
\begin{align}
  & \Phih := 
  \ui (\sh-s)(\ih-i) 
  + \hi (\sh-s) \ih
  + \hi s (\ih-i)
  \non
  \\
  & \quad {}
  + \ue (\sh-s)(\eh-e) 
  + \he (\sh-s) \eh
  + \he s (\eh-e) \,.
  \non
\end{align}
Then, a boring but trivial calculation shows that $\soluzh$ satisfies the variational equations
\begin{align}
  & \< \dt\xih , v >
  + \iO \bigl(
    (\ui i)\xih 
    + (\ui s) \iotah 
    + (\ue e)\xih 
    + (\ue s) \etah 
  \bigr) \, v
  \non
  \\
  & \quad {}
  + \iO \kas \nabla\xih \cdot \nabla v
  - \iO \gamma \rhoh \, v
  = - \iO \Phih \, v
  \non
  \\
  & \< \dt\etah , v >
  - \iO \bigl(
    (\ui i)\xih 
    + (\ui s) \iotah 
    + (\ue e)\xih 
    + (\ue s) \etah 
  \bigr) \, v
  + \iO (\sigma + \phie) \etah \, v
  \non
  \\
  & \quad {}
  + \iO \kae \nabla\etah \cdot \nabla v
  = \iO \Phih \, v
  \non
  \\
  & \< \dt \iotah , v >
  + \iO \phir \iotah \, v
  + \iO \kai \nabla\iotah \cdot \nabla v
  - \iO \sigma \etah \, v
  = 0
  \non
  \\
  & \< \dt\rhoh , v >
  - \iO \bigl( \phir \iotah + \phie \etah \bigr) \, v 
  + \iO \kar \nabla\rhoh \cdot \nabla v
  + \iO \gamma \rhoh \, v
  = 0
  \non
\end{align}
\aet\ and for every $v\in V$, 
as well as the initial condition
\Beq
  \soluzh(0) = (0,0,0,0) \,.
  \non
\Eeq
Now, we test the above equations by $\xih$, $\etah$, $\iotah$ and~$\rhoh$, respectively,
sum up and integrate over~$(0,t)$.
Then, we rearrange in order that the \lhs\ becomes
(recall the definition \eqref{defQt} of~$Q_t$)
\begin{align}
  & \frac 12 \iO |\xih(t)|^2
  + \frac 12 \iO |\etah(t)|^2
  + \frac 12 \iO |\iotah(t)|^2
  + \frac 12 \iO |\rhoh(t)|^2
  \non
  \\
  & \quad {}
  + \intQt \kas |\nabla\xih|^2
  + \intQt \kae |\nabla\etah|^2
  + \intQt \kai |\nabla\iotah|^2
  + \intQt \kar |\nabla\rhoh|^2 .
  \non
\end{align}
We have to estimate the \rhs\ we obtain.
We notice once more that all of the functions $\gamma$, $\ui$, $\ue$, $s$, $e$, $i$, $r$, $\sh$, $\eh$, $\ih$ and $\rh$
are bounded by known constants, so that many of the terms 
can be dealt with by simply owing to the Young inequality.
The more delicate terms come from the addends of $\Phih$ and we just consider two of them
since the other are analogous.
By recalling the compact embedding $V\emb\Lx4$, 
and applying first the \Holder, Young, Sobolev and compacness inequalities 
(see~\eqref{sobolev} and~\eqref{compact})
and then the continuous dependence estimate \eqref{contdep} 
(where we can replace the norm on the \rhs\ by the $\calU$~norm) 
applied to the states corresponding to $u$ and~$u+h$
we have~that
\begin{align}
  & - \intQt \ui (\sh-s) (\ih-i) \xih 
  \leq \iot \norma\ui_\infty \norma{(\sh-s)(\tau)} \, \norma{(\ih-h)(\tau)}_4 \, \norma{\xih(\tau)}_4 \, d\tau
  \non
  \\
  & \leq \iot \norma{\xih(\tau)}_4^2 \, d\tau
  + c \iot \norma{(\sh-s)(\tau)}^2 \,  \normaV{(\ih-h)(\tau)}^2
  \non
  \\
  \separa
  & \leq \frac \kamin 2 \intQt |\nabla\xih|^2
  + c \intQt |\xih|^2  
  + c \, \norma{\sh-s}_{\L\infty H}^2 \norma{\ih-h}_{\L2V}^2
  \non
  \\
  & \leq \frac \kamin 2 \intQt |\nabla\xih|^2
  + c \intQt |\xih|^2  
  + c \, \norma h_\calU^4 \,.
  \non  
\end{align}
The next term we consider is simpler.
We have~that
\begin{align}
  & - \intQt \hi (\sh-s) \ih \, \xih
  \leq \intQt |\xih|^2
  + c \, \norma\hi_\infty^2 \norma{\sh-s}_{\L2H}^2
  \leq \intQt |\xih|^2
  + c \, \norma h_\calU^4 \,.
  \non
\end{align}
By analogously treating the other terms and applying the Gronwall lemma,
we conclude~that
\Beq
  \norma\soluzh_\calY
  \leq c \, \norma h_\calU^{\pier 2} \,.
  \non
\Eeq
Since this inequality implies \eqref{tesifrechet}, the proof is complete.
\Edim


\subsection{First-order optimality conditions}
\label{OPT}

We introduce the functionals
$\calJ_1:(\C0H)^4\to\erre$ and $\calJ_2:(\LQ\infty)^2\to\erre$
by setting
\begin{align}
  & \calJ_1\soluz := 
  := \frac 12 \iO \bigl( \bigl( (e+i)(T) - \Lam \bigr)^+ \bigr)^2
  \non
  \\
  & \calJ_2(\ui,\ue)
  := \frac 12 \intQ \bigl( |\ui|^2 + |\ue|^2 \bigr) 
  \non
\end{align}
and we notice that they are \frechet\ differentiable.
Since
\Beq
  \calJ\uie = \calJ_1(\calS\uie,\ui,\ue) + \calJ_2\uie
  \quad \hbox{for every $\uie\in\Upiu$}
  \non
\Eeq
and Theorem~\ref{Frechet} holds,
we can compute the derivative of $\calJ$ at points of $\Upiu$ 
by applying the chain rule to the first addend.
By also observing that $\Uad$ is convex,
we have the following result:

\Bprop
\label{Badopt}
Assume that $\uiestar$ and $\soluzstar:=\calS\uiestar$
are an optimal control and the corresponding optimal state.
Then, the variational inequality
\Beq
  \iO \bigl( (\estar+\istar)(T) - \Lam \bigr)^+ \, (\eta+\iota)(T)
  + \intQ (\uistar\hi + \uestar\he) 
  \geq 0
  \label{badopt}
\Eeq
holds true for every $\uie\in\Uad$,
where $\eta$ and $\iota$ are the component\pier{s} of the solution $\soluzl$ to the linearized system \Pbll\
corresponding to the control pair $\uiestar$ and \pier{to} the variation $\hie$ given by
$\hi=\ui-\uistar$ and $\he=\ue-\uestar$.
\Eprop

We notice that the variation $\hie$ is admissible in the computation of the derivative of $\calS$
since $\uiestar+\hie=\uie$ belongs to~$\Upiu$.
However, this result is not satisfactory since it involves the linearized system infinitely many times.
This inconvenient is bypassed by introducing a proper adjoint problem.

Given an optimal control $\uiestar$ and the corresponding state $\soluzstar$,
the associated adjoint problem consists in looking for a quadruplet $\soluza$
that satisfies \pier{the} regularity requirement
\Beq
  p ,\,, q ,\,, w ,\,, z \in \H1\Vp \cap \L2V \emb \C0H
  \label{regsoluza}
\Eeq
and the variational equations
\begin{align}
  & - \< \dt p , v >
  + \iO \kas \nabla p \cdot \nabla v
  + \iO \bigl( \uistar \istar + \uestar \estar \bigr) (p - q) \, v
  = 0
  \label{primaa}
  \\
  \separa
  & - \< \dt q , v >
  + \iO \kae \nabla q \cdot \nabla v
  + \iO \bigl( \sigma + \phie - \uestar \sstar \bigr) q \, v
  \non
  \\
  & \quad {}
  + \iO \bigl( \uestar \sstar p - \sigma w - \phie z \bigr) \, v
  = 0
  \label{secondaa}
  \\
  \separa
  & - \< \dt w , v >
  + \iO \kai \nabla w \cdot \nabla v
  + \iO \bigl( \phir (w - z) + \uistar \sstar (p - q) \bigr) \, v
  = 0
  \label{terzaa}
  \\
  & - \< \dt z , v >
  + \iO \kar \nabla z \cdot \nabla v
  + \iO \gamma (z - p) \, v
  = 0
  \label{quartaa}
\end{align}
\aet\ and for every $v\in V$,
as well as the final conditions
\Beq
  p(T) = z(T) = 0
  \aand
  q(T) = w(T) = \bigl( (\estar+\istar)(T) - \Lam \bigr)^+ \,.
  \label{cauchya}
\Eeq
\Accorpa\Pbla primaa cauchya
Since all the coefficient\pier{s} are bounded\pier{, the final conditions are in $H$} and the diffusion coefficients satisfy~\eqref{hpk}, 
this is a standard backward parabolic problem.
Hence, it has a unique solution $\soluza$ with the regularity \eqref{regsoluza}.

At this point, we are ready to proof a satisfactory necessary condition for optimality.
Our statement involves the \pier{closed} convex sets
\begin{align}
  & \calU_i := \graffe{v\in\LQ2:\ 0 \leq v \leq \uimax\ \aeQ}
  \non
  \\
  & \calU_e := \graffe{v\in\LQ2:\ 0 \leq v \leq \uemax\ \aeQ}.
  \non
\end{align}

\Bthm
\label{Opt}
Let $\uiestar$ and $\soluzstar$ be an optimal control and the corresponding state, respectively,
and let $\soluza$ be the solution to the associated adjoint problem.
Then, the variational inequality
\Beq
  \intQ \bigl( \sstar\istar(q-p) + \uistar \bigr) (\ui-\uistar)
  + \intQ \bigl( \sstar\estar(q-p) + \uestar \bigr) (\ue-\uestar)
  \geq 0
  \label{opt}
\Eeq
holds true for every $\uie\in\Uad$ and
$\uistar$ and $\uestar$ are the $L^2$ projections 
of $\sstar\istar(p-q)$ and $\sstar\estar(p-q)$ on $\calU_i$ and $\calU_e$, respectively.
\Ethm

\Bdim
We consider the linearized problem \Pbll\ 
where $\uie$ and $\soluz$ are read $\uiestar$ and $\soluzstar$, respectively,
and the variation $\hie$ is generic for a while,~i.e.,
\begin{align}
  & \< \dt\xi , v >
  + \iO \bigl(
    (\uistar\istar)\xi 
    + (\uistar\sstar) \iota 
    + (\uestar\estar)\xi 
    + (\uestar\sstar) \eta 
  \bigr) \, v
  \non
  \\
  & \quad {}
  + \iO \kas \nabla\xi \cdot \nabla v
  - \iO \gamma \rho \, v
  = - \iO \bigl( (\sstar\istar) \hi + (\sstar\estar) \he \bigr) \, v
  \non
  \\
  \separa
  & \< \dt\eta , v >
  - \iO \bigl(
    (\uistar\istar)\xi 
    + (\uistar\sstar) \iota 
    + (\uestar\estar)\xi 
    + (\uestar\sstar) \eta 
  \bigr) \, v
  + \iO (\sigma + \phie) \eta \, v
  \non
  \\
  & \quad {}
  + \iO \kae \nabla\eta \cdot \nabla v
  = \iO \bigl( (\sstar\istar) \hi + (\sstar\estar) \he \bigr) \, v
  \non
  \\
  \separa
  & \< \dt \iota , v >
  + \iO \phir \iota \, v
  + \iO \kai \nabla\iota \cdot \nabla v
  - \iO \sigma \eta \, v
  = 0
  \non
  \\
  & \< \dt\rho , v >
  - \iO \bigl( \phir \iota + \phie \eta \bigr) \, v 
  + \iO \kar \nabla\rho \cdot \nabla v
  + \iO \gamma \rho \, v
  = 0
  \non
  \\
  & \soluzl(0) = (0,0,0,0)
  \non
\end{align}
where the variational equations hold \aet\ and for every $v\in V$.
Now, we fix an arbitrary element $\uie\in\Uad$ and choose 
$\hi=\ui-\uistar$ and $\he=\ue-\uestar$ in the above system.
Then, we test the equations by $p$, $q$, $w$ and~$z$, respectively,
sum up, integrate over~$(0,T)$ and rearrange.
We obtain~that
\begin{align}
  & \ioT \bigl(
    \< \dt\xi , p >
    + \< \dt\eta , q >
    + \< \dt \iota , w >
    + \< \dt\rho , z >
  \bigr) \, dt
  \non
  \\
  & \quad {}
  + \intQ \bigl(
    \kas \nabla\xi \cdot \nabla p
    + \kae \nabla\eta \cdot \nabla q
    + \kai \nabla\iota \cdot \nabla w
    + \kar \nabla\rho \cdot \nabla z
  \bigr)
  \non
  \\
  & \quad {}
  + \intQ \bigl(
    (\uistar\istar)\xi 
    + (\uistar\sstar) \iota 
    + (\uestar\estar)\xi
  \bigr) \, p 
  + \intQ (\uestar\sstar) \eta \, p
  - \intQ \gamma \rho \, p
  \non
  \\
  & \quad {}
  - \intQ \bigl(
    (\ui i)\xi 
    + (\uistar\sstar) \iota 
    + (\uestar\estar)\xi 
  \bigr) \, q
   \pier{{}+ \intQ (\uestar\sstar) \eta  \, q{}}
  + \intQ (\sigma + \phie) \eta \, q
  \non
  \\
  & \quad {}
  + \intQ \phir \iota \, w
  - \intQ \sigma \eta \, w
  - \intQ \bigl( \phir \iota + \phie \eta \bigr) \, z
  + \intQ \gamma \rho \, z
  \non
  \\
  & = - \intQ \bigl( (\sstar\istar) (\ui-\uistar) + (\sstar\estar) (\ue-\uestar) \bigr) \, (p - q) \,.
  \label{ltested}
\end{align}
At the same time, we consider the adjoint system \Pbla\ and test the equations 
by $-\xi$, $-\eta$, $-\iota$ and~$-\rho$, respectively.
Then, we add the equality we obtain to each other, integrate over~$(0,T)$ and rearrange.
We have~that
\begin{align}
  & \ioT \bigl(
    \< \dt p , \xi >
    + \< \dt q , \eta >
    + \< \dt w , \iota >
    + \< \dt z , \rho >
  \bigr) \, dt
  \non
  \\
& \quad {}  
  - \intQ \bigl(
    \kas \nabla p \cdot \nabla\xi
    + \kae \nabla q \cdot \nabla\eta
    + \kai \nabla w \cdot \nabla\iota
    + \kar \nabla z \cdot \nabla\rho
  \bigr)
  \non
  \\
  & \quad {}
  - \intQ \bigl( \uistar \istar + \uestar \estar \bigr) (p - q) \, \xi
  - \intQ \bigl( \sigma + \phie - \uestar \sstar \bigr) q \, \eta
  - \intQ \bigl( \uestar \sstar p - \sigma w - \phie z \bigr) \, \eta
  \non
  \\
  & \quad {}
  - \intQ \bigl( \phir (w - z) + \uistar \sstar (p - q) \bigr) \, \iota
  - \intQ \gamma (z - p) \, \rho
  = 0 \,.
  \label{atested}
\end{align}
At this point, we take the sum of \eqref{ltested} and \eqref{atested} 
and notice \pier{that} several cancellations \pier{occur}.
Namely, we obtain that
\begin{align}
  & \ioT \bigl(
    \< \dt p , \xi >
    + \< \dt q , \eta >
    + \< \dt w , \iota >
    + \< \dt z , \rho >
  \bigr) \, dt
  \non
  \\
  & \quad {}
  + \ioT \bigl(
    \< \dt p , \xi >
    + \< \dt q , \eta >
    + \< \dt w , \iota >
    + \< \dt z , \rho >
  \bigr) \, dt
  \non
  \\
  & = - \intQ \bigl( (\sstar\istar) (\ui-\uistar) + (\sstar\estar) (\ue-\uestar) \bigr) \, (p - q) \,.
  \non
\end{align}
By applying the \pier{well-known} integration-by-parts formula for functions belonging to the space $\H1\Vp\cap\L2V$
and accounting for the initial conditions for $\soluzl$ and the final conditions for $\soluza$,
we deduce~that
\Beq
  \iO \bigl( (\estar+\istar)(T) - \Lam \bigr)^+ (\eta+\iota)(T)
  = - \intQ \bigl( (\sstar\istar) (\ui-\uistar) + (\sstar\estar) (\ue-\uestar) \bigr) \, (p - q) \,.
  \non
\Eeq
By using this equality in \eqref{badopt} we obtain~\eqref{opt}.
To prove the last sentence of the statement, we observe that $\Uad=\calU_i\times\calU_e$.
Therefore, \eqref{opt} is equivalent to the pair of conditions
\begin{align}
  & \intQ \bigl( \sstar\istar(q-p) + \uistar \bigr) (v-\uistar) \geq 0
  \quad \hbox{for every $v\in\calU_i$}
  \non
  \\
  & \intQ \bigl( \sstar\estar(q-p) + \uestar \bigr) (v-\uestar) \geq 0
  \quad \hbox{for every $v\in\calU_e$}
  \non
\end{align}
and this concludes the proof of Theorem~\ref{Opt}. 
\Edim

\section*{Acknowledgments}
\pier{This research activity has been performed in the framework of the Italian-Romanian
collaboration agreement \textquotedblleft{Analysis and control of mathematical models 
for the evolution of epidemics, tumors and phase field processes}\textquotedblright\ between the
Italian CNR and the Romanian Academy.}
In addition, PC and ER gratefully mention some other support 
from the MIUR-PRIN Grant 2020F3NCPX 
``Mathematics for industry 4.0 (Math4I4)'' and the GNAMPA (Gruppo Nazionale per l'Analisi Matematica, 
la Probabilit\`a e le loro Applicazioni) of INdAM (Isti\-tuto 
Nazionale di Alta Matematica).  \pier{GM acknowledges the support of a grant of the Ministry of Research,
Innovation and Digitization, CNCS - UEFISCDI, project number
PN-III-P4-PCE-2021-0921, within PNCDI III.}



\Begin{thebibliography}{10}

\betti{%
\bibitem{acgrr} 
F. Auricchio, P. Colli, G. Gilardi, A. Reali, E. Rocca,
Well-posedness for a diffusion-reaction compartmental model simulating the spread of COVID-19,
{\it  Math. Models Methods Appl. Sci.} (2023) 1-20, DOI: 10.1002/mma.9196.}

\pier{%
\bibitem{Balderrama} R. Balderrama, J. Peressutti, J.P. Pinasco, F. Vazquez, C. S\'anchez de la Vega, 
Optimal control for a SIR epidemic model with limited quarantine, {\it Scientifc Reports} 
{\bf 12} 12583 (2022). DOI: 10.1038/s41598-022-16619-z.}

\betti{%
\bibitem{Bellomo4}
N. Bellomo, F. Brezzi, M.A.J. Chaplain,
Special Issue on ``Mathematics Towards COVID19 and Pandemic'',
{\it Math. Mod. Meth. Appl. Sci.} {\bf 31} (2021) issue 12.}

\betti{%
\bibitem{CGLMRR} 
P. Colli, H. Gomez, G. Lorenzo, G. Marinoschi, A. Reali, E. Rocca, 
Optimal control of cytotoxic and antiangiogenic therapies on prostate cancer growth., 
{\it Math. Models Methods Appl. Sci.} {\bf 31} (2021) 1419-1468.}

\pier{%
\bibitem{Mottoni-Orlandi-Tesei} P. de Mottoni, E. Orlandi, A. Tesei,
Asymptotic behavior for a system describing epidemics with migration and
spatial spread of infection, {\it Nonlinear Anal.}  {\bf 3}  (1979) 663-675.}

\pier{%
\bibitem{Fitzgibbon-1} W.E. Fitzgibbon, M.E. Parrott, G.F. Webb, Diffusion
epidemic models with incubation and crisscross dynamics, {\it Math.
Biosci.} {\bf 128} (1995) 131-155.}

\pier{%
\bibitem{Fitzgibbon-2} W.E. Fitzgibbon, J.J. Morgan, G.F. Webb, Y. Wu, A
vector-host epidemic model with spatial structure and age of infection,
{\it Nonlinear Anal. Real World Appl.} {\bf 41} (2018) 692-705.}

\betti{%
\bibitem{Gatto}
M. Gatto, E. Bertuzzo, L. Mari, S. Miccoli, L. Carraro, R. Casagrandi, A. Rinaldo,
Spread and dynamics of the COVID-19 epidemic in Italy: Effects of emergency containment measures,
{\it Proc. Nat. Acad. Sci.} {\bf 117} (2020) 10484-10491.}

\betti{%
\bibitem{Giordano}
G. Giordano, F. Blanchini, R. Bruno, P. Colaneri, A. Di Filippo, A. Di Matteo, M. Colaneri,
Modelling the COVID-19 epidemic and implementation of population-wide interventions in Italy,
{\it Nat. Med.} {\bf 26} (2020) 855-860.}

\betti{%
\bibitem{Jha}
P.K. Jha, L. Cao, J.T. Oden, 
Bayesian-based predictions of COVID-19 evolution in Texas using multispecies mixture-theoretic continuum models,
{\it Comput. Mech.} {\bf 66} (2020) 1055-1068.}

\bibitem{LSU}
O.\,A. Lady\v{z}enskaja, V.\,A. Solonnikov and N.\,N. Uralceva,
``Linear and Quasilinear Equations of Parabolic
Type'', Mathematical Monographs, Vol.~{\bf 23}, 
American Mathematical Society, Providence, Rhode Island, 1968.

\betti{%
\bibitem{Linka}
K. Linka, P. Rahman, A. Goriely, E. Kuhl,
Is it safe to lift COVID-19 travel bans? The Newfoundland story,
{\it Comput. Mech.} {\bf 66} (2020) 1081-1092.}

\pier{%
\bibitem{GM-AMO} G. Marinoschi, Parameter estimation of an epidemic model
with state constraints, {\it Appl. Math. Optimiz.} {\bf 84} (2021) suppl. 2, S190-S1923.}

\pier{%
\bibitem{GM-DCDS} G. Marinoschi, Identification of transmission rates and
reproduction number in a SARS-CoV-2 epidemic model, {\it Discrete Contin. Dyn.
Syst. Ser. S} {\bf 15} (2022) 3735-3744.}

\pier{%
\bibitem{Medhaoui} M. Mehdaoui, A.L. Alaoui, M. Tilioua,  Optimal control
for a multi-group reaction-diffusion SIR model with heterogeneous incidence
rates, {\it Int. J. Dynam. Control} (2022). DOI: 10.1007/s40435-022-01030-3.}

\bibitem{Simon}
J. Simon,
{Compact sets in the space $L^p(0,T; B)$},
{\it Ann. Mat. Pura Appl.~(4)}
{\bf 146} (1987), 65-96.

\betti{%
\bibitem{Vig1}
A. Viguerie, G. Lorenzo, F. Auricchio, D. Baroli,
T.J.R. Hughes, A. Patton, A. Reali,
Th.E. Yankeelov, A. Veneziani,
Simulating the spread of COVID-19 via a spatially-resolved
susceptible-exposed-infected-recovered-deceased (SEIRD) model with heterogeneous diffusion,
{\it Appl. Math. Lett.} {\bf 111} (2021) Paper No. 106617, 9 pp.
}

\betti{%
\bibitem{Vig2}
A. Viguerie, A. Veneziani, G. Lorenzo, D. Baroli, N. Aretz-Nellesen,
A. Patton, T.E. Yankeelov, A. Reali, T.J.R. Hughes, F. Auricchio,
Diffusion-reaction compartmental models formulated in a continuum mechanics framework: application to COVID-19, mathematical analysis, and numerical study,
{\it Comput. Mech.} {\bf 66} (2020) 1131-1152.}

\betti{%
\bibitem{Wang}
Z. Wang, X. Zhang, G.H. Teichert, M. Carrasco-Teja, K. Garikipati,
System inference for the spatio-temporal evolution of infectious diseases: Michigan in the time of COVID-19,
{\it Comput. Mech.} {\bf 66} (2020) 1153-1176.}

\pier{%
\bibitem{Webb} G.F. Webb, A reaction-diffusion model for a deterministic
diffusive epidemic, {\it J. Math. Anal. Appl.} {\bf 84} (1981) 150-161.}

\betti{%
\bibitem{Zohdi}
T.I. Zohdi,
An agent-based computational framework for simulation of global pandemic and social response on planet X,
{\it Comput. Mech.} {\bf 66} (2020) 1195-1209.}

\End{thebibliography}

\End{document}
